\documentclass[12pt]{article}
\usepackage[left=108pt, right=74pt, top=108pt, bottom=81pt, dvips, pdftex]{geometry}
\usepackage{amsmath, amssymb, amsthm, amsfonts,graphicx}
\usepackage{epstopdf}

\newtheorem{theorem}{Theorem}
\newtheorem{lemma}{Lemma}

\newtheorem{proposition}{Proposition}
\theoremstyle{definition}

\def\ad{\mathop{\rm ad}}
\def\Ker{\mathop{\rm Ker}}
\def\g{\mathfrak g}
\def\Oplus{\bigoplus\limits}
\title{Laplacian spectrum for the nilpotent Kac-Moody Lie algebras}
 \author{Dmitry Fuchs \\University of California\\Davis
\and Constance Wilmarth\\University of California\\Davis}
 \date{}

\begin{document}
\maketitle
\begin{abstract}
\noindent
We prove that the maximal nilpotent subalgebra of a Kac-Moody Lie algebra has a (essentially, unique) Euclidean metric with respect to which the Laplace operator in the chain complex is scalar on each component of a given degree. Moreover, both the Lie algebra structure and the metric are uniquely determined by this property.
\end{abstract}
\section{Introduction}
Let $\g$ be a real Lie algebra which is either finite-dimensional or has a grading $\g=\Oplus_{{\bf k}\in{\mathbb Z}^n}\g^{({\bf k})}$ such that all the chain spaces $C_q^{({\bf k})}(\g)=\bigoplus\limits_{{\bf k}_1+\dots+{\bf k}_q={\bf k}}(\g^{({\bf k}_1)}\wedge\dots\wedge\g^{({\bf k}_q)})$ are finite-dimensional. (Below we consider only the case when $\g=\bigoplus\limits_{(k_1,\dots,k_n)\succ(0,\dots,0)}\g^{(k_1,\dots,k_n)}$ where the notation $(k_1,\dots,k_n)\succ(\ell_1,\dots,\ell_n)$ means that $k_1\ge\ell_1,\dots,k_n\ge\ell_n,$ and $(k_1,\dots,k_n)\ne(\ell_1,\dots,\ell_n)$, and all the spaces $\g^{(k_1,\dots,k_n)}$ are finite-dimensional.) Suppose that for each value of $\bf k$, some Euclidean structure is fixed for $\g^{(\bf k)}$. Then Euclidean structures arise in all the chain spaces $C_q^{(\bf k)}(\g)$, and they give rise to canonical isomorphisms between the chain spaces $C_q^{(\bf k)}(\g)$ and the corresponding cochain spaces, $C^q_{(\bf k)}(\g)=(C_q^{(\bf k)}(\g))^\ast$. Thus, we can regard the boundary and coboundary operators as acting in the same spaces, $\partial\colon C_q^{(\bf k)}(\g)\to C_{q-1}^{(\bf k)}(\g),\ \delta\colon C_q^{(\bf k)}(\g)\to C_{q+1}^{(\bf k)}(\g)$, and to form the {\it Laplace operators} $\Delta\colon C_q^{(\bf k)}(\g)\to  C_q^{(\bf k)}(\g)$. Chains (cochains) annihilated by $\Delta$ are called {\it harmonic}. The finite-dimensional version of the Hodge--de Rham theory yields the following result.

\begin{proposition}

Every harmonic chain (cochain) is a cycle (cocycle), and every homology (cohomology) class of $\g$ (with trivial coefficients) is represented by a unique harmonic chain (cochain). In particular, there are canonical isomorphisms $$\Ker[\Delta\colon C_q^{(\bf k)}(\g)\to  C_q^{(\bf k)}(\g)]=H_q^{(\bf k)}(\g)=H^q_{(\bf k)}(\g).$$

\end{proposition}

(For details, see \cite{Fuchs}, Section 1.5.3.)\smallskip

\paragraph{Remark.} The results discussed below indicate that not only the kernel, but the whole spectrum of the Laplacian must have a significance for the (co)homology. However, this significance is not clear to us.

To our best knowledge, the spectrum of the Laplacian has been calculated in two cases. 

First, it is known for the Lie algebra $L_1(1)$ of polynomial vector fields in the line with an at least double zero at 0. This algebra has a basis $\{e_i\mid i>0\}$ with the commutator operation $[e_i,e_j]=(j-i)e_{i+j}$. We introduce in this algebra a $\mathbb Z$-grading and a Euclidean structure letting $\deg e_i=i,\ \|e_i\|=1$. For positive integers $i_1,\dots,i_q$ such that $i_r-i_{r-1}\ge3$ for $r=2,\dots,q$, let \[ \begin{array} {rl} E(i_1,\dots,i_q)&=\displaystyle{\sum_{s=1}^q{i_s\choose3}-\sum_{1\le \ell<m\le q}i_\ell i_m,}\\ \alpha_r(i_1,\dots,i_q)&=\left\{\displaystyle{\begin{array} {rl} 0,&\mbox{if}\ r=1,i_1<3,\\ 1,&\mbox{if}\ r=1,i_1\ge3,\\ 0,&\mbox{if}\ 1<r\le q,i_r-i_{r-1}=3,\\ 1,&\mbox{if}\ 1<r\le q,i_r-i_{r-1}>3 \end{array} }\right.\\\alpha(i_1,\dots,i_q)&=\displaystyle{\sum_{r=1}^q\alpha_r(i_1,\dots,i_q)}\end{array} \] It is easy to check that $E(1,4,7,\dots,3q-2)=E(2,5,8,\dots,3q-1)=0$, and all other values of the function $E$ are positive.

\begin{theorem} [\cite{GFF}, \cite{Wein}]

The set of eigenvalues of the Laplace operator $\Delta\colon C_\ast(L_1(1))\to C_\ast(L_1(1))$ coincides with the set of numbers $E(i_1,\dots,i_q)$. The multiplicity of the eigenvalue $E(i_1,\dots,i_q)$ equals $2^{\alpha(i_1,\dots,i_q)}$. (Occasional coincidences $E(i_1,\dots,i_q)=E(i'_1,\dots,i'_{q'})$ are possible; in such cases the multiplicities are added.)

\end{theorem}

(For a sketch of a proof see \cite{Fuchs}, Section 2.3.1(B).)\smallskip

Second, it is known for the nilpotent current algebra $T^+(n)$ which is described in the following way. Consider a Lie algebra $\widetilde T^+(n)$ of $(n\times n)$-matrices $\|p_{ij}(t)\|,\ p_{ij}(t)\in{\mathbb R}[t]$ such that $p_{ij}(0)=0$, if $i\ge j$, with the usual commutator operation. Let $E_{ij}^r\in\widetilde T^+(n)$ be a matrix with the only non-zero entry $p_{ij}(t)=t^r$ (thus, $1\le i\le n,1\le j\le n, r\ge0$ and $r>0$ ,if $i\ge j$). The Lie algebra $\widetilde T^+(n)$ has a natural $n$-grading, \[\deg E_{ij}^r=
\left\{ \begin{array} {rl} (\underbrace{r,\dots,r}_{i-1},\underbrace{r+1,\dots,r+1}_{j-i},\underbrace{r,\dots,r}_{n-j+1}),&\mbox{if}\ i\le j,\\ (\underbrace{r,\dots,r}_{j-1},\underbrace{r-1,\dots,r-1}_{i-j},\underbrace{r,\dots,r}_{n-i+1}),&\mbox{if}\ i>j, \end{array} \right.\] and a natural Euclidean structure for which $\{E_{ij}^r\}$ is an orthonormal basis. We set $$T^+(n)=\left\{\|p_{ij}\|\in\widetilde T^+(n)\mid p_{11}(t)+\dots+p_{nn}(t)=0\right\}\subset\widetilde T^+(n);$$ obviously, $T^+(n)$ inherits from $\widetilde T^+(n)$ the structure of a Lie algebra, the grading, and the Euclidean structure.

\begin{theorem}[\cite{Fei},\cite{Le}]

The Laplace operator $\Delta\colon C_\ast^{(k_1,\dots,k_n)}(T^+(n))\to C_\ast^{(k_1,\dots,k_n)}(T^+(n))$ is the multiplication by $-\sum_ik_i^2+\sum_ik_ik_{i+1}+\sum_ik_i$ (we set $k_{n+1}=k_1$).

\end{theorem}

(For a sketch of a proof see \cite{Fuchs}, Section 2.5.1.)\smallskip

The goal of this paper is to provide a generalization of Theorem 2 to the case of the maximal nilpotent subalgebra of an arbitrary Kac-Moody algebra. (It should be mentioned that no generalization, or explanation, exists for Theorem 1.) We supply below all the necessary definitions; for the general theory of Kac-Moody Lie algebras see \cite{Kac}.

Let $A=\|a_{ij}\|$ be an $n\times n$ matrix with all the diagonal entries equal to 2 and all non-diagonal entries being non-positive integers. We assume the matrix $A$ {\it symmetrizable} which means that there exist a diagonal matrix $D$ whose diagonal entries $d_1,\dots,d_n$ are positive integers such that the matrix $DA$ is symmetric. We may also assume the matrix $A$ {\it irreducible} which means that there is no partition of $\{1,\dots,n\}$ into non-empty parts $I,J$ such that $a_{ij}=0$ for all $i\in I,j\in J$. Let $G=G(A)$ be the (real) Kac-Moody Lie algebra with the Cartan matrix $A$, and let $N=N(A)$ be the corresponding nilpotent Lie algebra. In other words, $N$ has a system of generators $e_1,\dots,e_n$ with the defining set of relations $(\ad e_i)^{-a_{ij}+1}e_j=0$. The algebra $N$ has a natural $n$-grading, $N=\Oplus_{(k_1,\dots,k_n)\succ(0,\dots,0)}N^{(k_1,\dots,k_n)}$ where $N^{(k_1,\dots,k_n)}$ consists of linear combinations of commutator monomials of the generators involving precisely $k_i$ letters $e_i\ (i=1,\dots,n)$. The following statement is our main result.

\begin{theorem}

There exist unique Euclidean structures in the spaces $N^{(k_1,\dots,k_n)}$ such that $\|e_i\|=1\ (i=1,\dots,n)$ and the corresponding Laplace operator $\Delta\colon C^{(k_1,\dots,k_n)}_\ast(N)\to
C^{(k_1,\dots,k_n)}_\ast(N)$ is the multiplication by $$E(k_1,\dots,k_n)=\sum_id_ik_i-\frac12\sum_{i,j}d_ia_{ij}k_ik_j.$$

\end{theorem}

In the case when \[A=\left[ \begin{array} {cccccc} \phantom{-}2&-1&&&&-1\\ -1&\phantom{-}2&-1&&&\\ &-1&\phantom{-}2&\ddots&& \\ &&\ddots&\ddots&\ddots&\\ &&&\ddots&\phantom{-}2&-1\\ -1&&&&-1&\phantom{-}2\end{array}\right], \] this is equivalent to Theorem 2.

\begin{proposition}

If $E(k)\ne0$, then $H_\ast^{(k)}(N(A))=0$.

\end{proposition}

This follows from Proposition 1 and Theorem 3. \smallskip

Proposition 2 is not new: it is essentially contained in \cite{KK}. More precisely, \cite{KK} yields a description of a Bernstein-Gelfand-Gelfand resolution of the trivial module over a Kac-Moody Lie algebra. This is also a free resolution of the trivial module over $N(A)$.

\section{Proof of Main Theorem}

\subsection{The Laplace operator has order 2}

We begin by recalling the notion of the order of a differential operator in the standard calculus. A linear operator $D:C^\infty({\mathbb R})\to C^\infty({\mathbb R})$ is a differential operator of degree 1 (that is, $D(f)=af'+bf$ where $a$ and $b$ are functions), if the identity $$D(fg)=D(f)g+D(g)f-D(1)fg$$ holds for any functions $f,g$. Similarly, an operator of degree 2 is characterized by the identity $$D(fgh)=D(fg)h+D(fh)g+D(hg)f-D(f)gh-D(g)fh-D(h)fg+D(1)fgh$$ (and so on; but we do not need operators of orders greater than 2). It is well known that the commutator of operators of order $p$ and $q$ has the order $p+q-1$.

In the non-commutative (super-commutative) case of chains/cochains of a Lie algebra (with a Euclidean structure) , the notion of a differential order looks slightly different. In particular, the operator $\delta\colon C_\ast(\g)\to C_\ast(\g)$ has order 1 which means that$$\delta(c_1\wedge c_2)=\delta(c_1)\wedge c_2+(-1)^{d_1d_2}\delta(c_2)\wedge c_1$$ for $c_i\in C_{d_i}(\g)$. However, the operator $\partial\colon C_\ast(\g)\to C_\ast(\g)$ has order 2 which means that \[ \begin{array} {rl} &\partial(c_1\wedge c_2\wedge c_3)\\ &\hskip.4in =\partial(c_1\wedge c_2)\wedge c_3+(-1)^{d_2d_3}\partial(c_1\wedge c_3)\wedge c_2+(-1)^{d_1(d_2+d_3)}\partial(c_2\wedge c_3)\wedge c_1\\ &
\hskip.4in -\partial(c_1)\wedge c_2\wedge c_3-(-1)^{d_1d_2}\partial(c_2)\wedge c_1\wedge c_3-(-1)^{(d_1+d_2)d_3}\partial(c_3)\wedge c_1\wedge c_2 \end{array} \] for $c_i\in C_{d_i}(\g)$. Since the Laplace operator is a (super)commutator of $\partial$ and $\delta$, it also has order 2, and we have the following lemma.

\begin{lemma}

The Laplace operator $\Delta\colon C_\ast(\g)\to C_\ast(\g)$ has the order 2, that is, \[ \begin{array} {rl} &\Delta(c_1\wedge c_2\wedge c_3)\\ &\hskip.4in =\Delta(c_1\wedge c_2)\wedge c_3+(-1)^{d_2d_3}\Delta(c_1\wedge c_3)\wedge c_2+(-1)^{d_1(d_2+d_3)}\Delta(c_2\wedge c_3)\wedge c_1\\ &
\hskip.4in -\Delta(c_1)\wedge c_2\wedge c_3-(-1)^{d_1d_2}\Delta(c_2)\wedge c_1\wedge c_3-(-1)^{(d_1+d_2)d_3}\Delta(c_3)\wedge c_1\wedge c_2 \end{array} \] for all $c_i\in C_{d_i}(\g)$.
 
\end{lemma}

\paragraph{Remark.} It is important that Lemma 1 is compatible with Theorem 3 in the following sense: if 
$\g=N=N(A)$ and $c_i\in C_{d_i}^{(p_i)}(N)$ where $(p_i)=(p_{i1},\dots,p_{in})$, then every term in the equality of Lemma 1 is $c_1\wedge c_2\wedge c_3$ times an approrpiate eigenvalue of $\Delta$, and the equality becomes $$E(p_1+p_2+p_3)=E(p_1+p_2)+E(p_1+p_3)+E(p_2+p_3)-E(p_1)-E(p_2)-E(p_3)$$which is true (because $E$ is a polynomial of degree 2).

\subsection{Construction of a Lie algebra with a given Laplace operator}

Let $A, a_{ij}, D, d_i$ denote the same as in Section 1. We will now construct a graded Lie algebra $\g=\Oplus_{(k_1,\dots,k_n)\succ(0,\dots,0)}\g^{(k_1,\dots,k_n)}$ with Euclidean structures in (finite-dimensional) spaces $\g^{(k_1,\dots,k_n)}$ satisfying the conclusion of Theorem 3 (with $N$ replaced by $\g$). Moreover, we will see that $\g$ is unique up to an isometric isomorphism, provided that $\dim\g^{(1,0,\dots,0)}=\dim\g^{(0,1,0,\dots,0)}=\dots=\dim\g^{(0,\dots,0,1)}=1$. (Later on, we will see that $\g$=$N(A)$.)

First consider a given graded Lie algebra $\g$ with Euclidean structures in $\g^{(k_1,\dots,k_n)}$. Choose an orthonormal basis in each $\g^{(k)},\ (k)=(k_1,\dots,k_n)$; then wedge products of the elements of the bases in $\g^{(k)}$ form orthonormal bases in the chain spaces $(\Lambda^q\g)^{(k)}$. For a fixed $(k)=(k_1,\dots,k_n)\succ(0,\dots,0)$, consider the matrix

\begin{figure}[hbtp]
\centering
\includegraphics[width=3in]{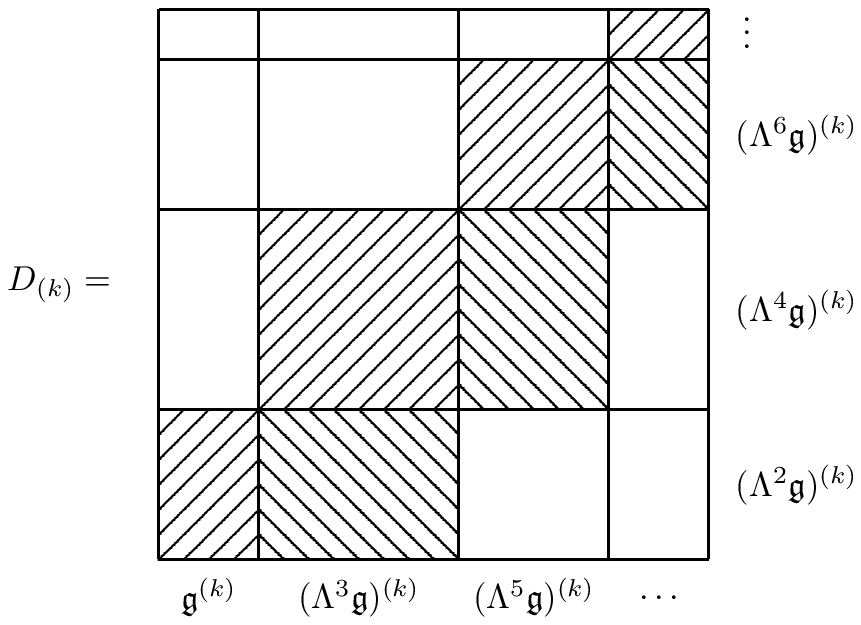}
\end{figure}

\noindent with rows (columns) labeled by the elements of our orthonormal bases in $(\Lambda^q\g)^{(k)}$ with $q$ even (odd). Let the shadowed blocks represent the boundary/coboundary operators in the chain/cochain complexes of $\g$, and let the unshadowed blocks be zero. Take two rows or two columns of the matrix $D_{(k)}$ corresponding to bases elements $c\in(\Lambda^q\g)^{(k)},\ c'\in(\Lambda^{q'}\g)^{(k)}$ (so $q$ and $q'$ have the same parity) and compute their dot-product. If $|q'-q|>2$, then this dot-product is obviously zero. If $|q'-q|=2$, it is also zero because of the relations $\partial\circ\partial=0,\ \delta\circ\delta=0$. Finally, if $q'=q$, then this dot-product is the coefficient at $c'$ in $\Delta(c)$ (and the coefficient at $c$ in $\Delta(c ')$).

Now, if the Laplace operator $\Delta\colon C_\ast^{(k)}\to C_\ast^{(k)}$ is the multiplication by a positive number $\lambda$, then the dot-product of every two different rows, as well as of every two different columns, is equal to zero, and the dot-square of every row or column is equal to $\lambda$; in other words, the whole matrix $D_{(k)}$ is an orthogonal matrix times $\sqrt\lambda$. 

This paves the way for a construction announced in the beginning of the section. First, we put $\dim\g^{(0,\dots,0,1,0,\dots,0)}=1$ and choose (in an arbitrary way) non-zero vectors $e_1\in\g^{(1,0,\dots,0)},\dots,e_n\in\g^{(0,\dots,0,1)}$ to have the length 1. Take a $(k)=(k_1,\dots,k_n)$ where $k_i$'s are non-negative integers with $k_1+\dots+k_n>1$. If $E(k)\le0$, we put $\g^{(k)}=0$; let $E(k)>0$. The matrix $D_{(k)}$ described above is fully determined, except the left bottom (shadowed) block. Away from this block, the dot product of every two distinct rows or columns is zero, and the dot-square of every row or column is equal to $E(k)$.  This follows from the identities $\partial\circ\partial=0,\, \delta\circ\delta=0$ and also from Lemma 1 and the remark after it which implies that the Laplace operator $\Delta\colon C^{(k)}_q\to C^{(k)}_q$ with $q\ge3$ (fully determined) is multiplication by $E(k)$. Thus, the columns of our matrix disjoint from the left bottom box are pairwise orthogonal and have dot squares $E(k)$. We can construct the missing columns making the whole matrix an orthogonal matrix times $\sqrt{E(k)}$. Since the dot-squares of the rows above the left bottom block are already equal to $E(k)$, the new columns will be confined to this block. Thus, we will have a $\g^{(k)}$ (with $\dim\g^{(k)}=\sum\limits_{q\ge2,\ {\rm even}}\dim(\Lambda^q\g)^{(k)}-\sum\limits_{q\ge3,\ {\rm odd}}\dim(\Lambda^q\g)^{(k)}$) with a ready orthonormal basis, and the new box yields a bracket $[\, ,\, ]\colon(\Lambda^2\g)^{(k)}\to\g^{(k)}$. Moreover, the orthogonality of the columns of the new box to the columns of the box next to the right means precisely that this bracket satisfies the Jacobi identity. (Notice that it could happen that $\sum\limits_{q\ge2,\ {\rm even}}\dim(\Lambda^q\g)^{(k)}=\sum\limits_{q\ge3,\ {\rm odd}}\dim(\Lambda^q\g)^{(k)}$; in this case we do not need any new columns, and simply put $\g^{(k)}=0$.)

This completes the construction promised in the beginning of the section; the uniqueness is obvious.

\subsection{End of the proof}

It remains to prove that the Lie algebra $\g$ of Section 2.2 is $N(A)$. This follows from three remarks.

First, it follows from the construction of Section 2.2 that if $(k_1,\dots,k_n)\succ(0,\dots,0)$ and $k_1+\dots+k_n\ge1$, then the bracket mapping $[\, ,\, ]\colon(\Lambda^2\g)^{(k)}\to\g^{(k)}$ is onto; hence, $\g$ (like $N(A)$) is generated by $e_1,\dots,e_n$.

Second, the defining relations $(\ad e_i)^{-a_{ij}+1}e_j=0$ hold. Indeed, the degree $(k)=(k_1,\dots,k_n)$ of $(\ad e_i)^{-a_{ij}+1}e_j$ is described by the equalities $k_i=-a_{ij}+1, k_j=1, k_s=0$ for $s\ne i,j$. Hence, \[ \begin{array} {rl} E(k)&=\displaystyle{\sum d_ik_i-\frac12 \sum a_{ij}k_ik_j}\\ &=d_i(-a_{ij}+1)+d_j-d_i(-a_{ij}+1)^2-d_j-d_ia_{ij}(a_{ij}+1)\\ &=-d_ia_{ij}+d_i+d_j-d_ia_{ij}^2+2d_ia_{ij}-d_i-d_j+d_ia_{ij}^2-d_ia_{ij}=0\end{array}\] By construction, this means that $\g^{(k)}=0$, hence $(\ad e_i)^{-a_{ij}+1}e_j=0$. Thus, there is a graded epimorphism $N(A)\to\g$.

Third, it is true that for all $(k),\ \dim\g^{(k)}=\dim N(A)^{(k)}$. Indeed, for any $(k)$ with $E(k)\ne0$, the dimensions $\dim\g^{(k)}$ are determined inductively from the relation $\sum(-1)^q\dim(\Lambda^q\g)^{(k)}=0$. A similar relation, $\sum(-1)^q\dim(\Lambda^qN(A))^{(k)}=0$ (for the same values of $(k)$) follows from Proposition 2 and the Euler-Poincar\'e Lemma. In addition to that, $\dim\g^{(k)}=\dim N(A)^{(k)}=1$, if $(k)=(0,\dots,0,1,0,\dots,0)$, and $\g^{(k)}=N(A)^{(k)}=0$, if $(k)=(k_1,\dots,k_n)\succ(0), k_1+\dots+k_n>1$, and $E(k)\le0$. Hence, our epimorphism $N(A)\to\g$ is, actually, an isomorphism.

\section{Conclusion}

\subsection{Canonical basis in $N(A)$}

The construction of Section 2.2 shows that the maximal nilpotent subalgebra of a Kac-Moody Lie algebra has a canonical Euclidean metric. The metric depends on the choice of generators of length 1, but the commutator relations do not depend on anything. In some cases (like Theorem 2) this metric looks  usual, but sometimes, even in the finite-dimensional case, it is less obvious. For example, the maximal nilpotent subalgebra of the rank 2 exceptional Lie algebra $G_2$ has dimension 6. The Cartan matrix is $A=\displaystyle{\left[\begin{array} {rr} 2&-1\\ -3&2 \end{array}\right]}$. There is a basis $\{e_{0,1},e_{1,0},e_{1,1}, e_{1,2},e_{1,3},e_{2,3}\},\ \deg e_{i,j}=(i,j)$ in $N(A)$ with the commutator relations 
\[ [e_{0,1},e_{1,0}]=\sqrt3\, e_{1,1},\, [e_{0,1},e_{1,1}]=2\, e_{1,2},\, [e_{0,1},e_{1,2}]=
\sqrt3\, e_{1,3},\]\vskip-25pt \[[e_{1,0},e_{1,3}]=\sqrt3\, e_{2,3},\, [e_{1,1},e_{1,2}]=\sqrt3\, e_{2,3}. \] If we regard this basis as orthonormal, then the Laplace operator in $C_\ast^{(p,q)}$ is the multiplication by $3p+q-3p^2-q^2+3pq$.

A more interesting example is provided by the twisted affine Kac-Moody Lie algebra $A^{(2)}_2$ with the Cartan matrix $A=\displaystyle{\left[\begin{array} {rr} 2&-1\\ -4&2 \end{array}\right]}$. This Lie algebra (after factoring over the one-dimensional center) is embedded into the current Lie algebra ${\mathfrak sl}(3)\otimes{\mathbb R}[t,t^{-1}]$. It is well known that it has a basis $e_i$ such that $[e_i,e_j]=\alpha_{ij}e_{i+j}$ where the numbers $\alpha_{ij}$ depend only on $i,j\bmod8$ (see Kac's book \cite{Kac}, Exercise 8.16).

The basis given in the book of Kac is not precisely our canonical basis; to get the latter, we need to modify  it by some coefficients: \[ e_{8s}=\sqrt2\left[\begin{array} {ccc} t^{2s}&0&0\\ 0&0&0\\ 0&0&-t^{2s}\end{array} \right];\quad e_{8s+1}=2\left[\begin{array} {ccc} 0&t^{2s}&0\\ 0&0&t^{2s}\\ 0&0&0\end{array} \right];\quad e_{8s+2}=\left[\begin{array} {ccc} 0&0&0\\ 0&0&0\\ t^{2s+1}&0&0\end{array} \right];\] \[e_{8s+3}=\left[\begin{array} {ccc} 0&0&0\\ t^{2s+1}&0&0\\ 0&-t^{2s+1}&0\end{array} \right];\quad e_{8s+4}=\sqrt{\frac23}\left[\begin{array} {ccc} t^{2s+1}&0&0\\ 0&-2t^{2s+1}&0\\ 0&0&t^{2s+1}\end{array} \right];\] \[e_{8s+5}=2\left[\begin{array} {ccc} 0&t^{2s+1}&0\\ 0&0&-t^{2s+1}\\ 0&0&0\end{array} \right];\quad e_{8s+6}=4\left[\begin{array} {ccc} 0&0&t^{2s+1}\\ 0&0&0\\ 0&0&0\end{array} \right];\] \[e_{8s+7}=\left[\begin{array} {ccc} 0&0&0\\ t^{2s+2}&0&0\\ 0&t^{2s+2}&0\end{array} \right].\]

The commutator of the elements of this basis is given by the formula  $[e_{8s+i},e_{8s'+j}]=\alpha_{ij}e_{8(s+s')+(i+j)},\ 1\le i\le8, 1\le j\le8$ with the $8\times8$ matrix $\|\alpha_{ij}\|$ being \[ \left[\begin{array} {rrrrrrrr} 0&2&\sqrt6&-\sqrt6&-2&0&\sqrt2&-\sqrt2\\ -2&0&0&0&2&-\sqrt8&0&\sqrt8\\ -\sqrt6&0&0&\sqrt6&-\sqrt2&2&-2&\sqrt2\\ \sqrt6&0&-\sqrt6&0&\sqrt6&0&-\sqrt6&0\\ 2&-2&\sqrt2&-\sqrt6&0&0&\sqrt6&-\sqrt2\\ 0&\sqrt8&-2&0&0&0&2&-\sqrt8\\ -\sqrt2&0&2&\sqrt6&-\sqrt6&-2&0&\sqrt2\\ \sqrt2&-\sqrt8&-\sqrt2&0&\sqrt2&\sqrt8&-\sqrt2&0\end{array} \right].\]

The natural grading of the Lie algebra $A^{(2)}_2$ is given by the following rule: if $-1\le s\le6$, then       \[ \deg e_{8n+s}=\left\{ \begin{array} {ll} (4n+s,2n),&\mbox{if}\  s\le1,\\ (4n+s-2,2n+1),&\mbox{if}\  s>1.\end{array}\right.\] The Laplace operator $\Delta\colon C_\ast^{(p,q)}\to C_\ast^{(p,q)}$ with respect to the metric determined by the basis $\{e_i,i>0\}$ is the multiplication by $4p+q-4p^2-q^2+4pq$.

\subsection{Some remarks on the multiplicative structure in $H^\ast(N(A))$}

It follows from our results (and, actually, can be proved directly) that there is a basis in $H^\ast(N(A))$ represented by uniquely chosen monomial cochains (that is, products of elements of the basis in $C^1(N(A))=N(A)^\ast$ dual to our canonical basis). This gives rise to a description of the multiplication in $H^\ast(N(A))$, which, however, is not very explicit. Let us begin with a couple of simple remarks.

First, it follows from the description above that the multiplication in $H^\ast(N(A))$ is ``square-free": the square of any cohomology class is zero.

Second, every monomial cochain representing a non-zero element of $H^\ast(N(A))$ should contain at least one factor from $C^1_{(0,\dots,0,1,0,\dots,0)}(N(A))$; this implies that the cohomological length of $H^\ast(N(A))$ does not exceed the rank of $G(A)$, that is the size of $A$.

Third, in the finite-dimensional case, the multiplication in $H^\ast(N(A))$ satisfies the Poincar\'e duality: if a non-zero element $\alpha\in H^q(N(A))$ is represented by a monomial cochain $c_{i_1}\dots c_{i_q}$, then the complimentary monomial $c_{j_1}\dots c_{j_r},\ q+r=d=\dim N(A)$ also represents a non-zero cohomology class, $\beta\in H^r(N(A))$, and $\alpha\beta$ is a non-zero element in $H^d(N(A))\cong{\mathbb R}$. It follows from the preceding remark that in the finite-dimensional case of rank 2 there are no other non-zero products. (It seems likely that in the infinite-dimensional case of rank 2, the multiplication in $H^\ast(N(A))$ is trivial.)

Now, let us consider some examples. Let $N(A)={\mathfrak n}(n)$ be the Lie algebra of (strictly) upper triangular $n\times n$ matrices, associated to the Cartan matrix \[A=\left[\begin{array} {rrrrr} 2&-1&&&\\ -1&2&-1&&\\ &\ddots&\ddots&\ddots&\\ &&-1&2&-1\\ &&&-1&2 \end{array}\right] \] For this Lie algebra, $\dim{\mathfrak n}(n)=\displaystyle{\frac{n(n-1)}2}$ and $\dim H^\ast({\mathfrak n}(n))=n!$. The basis in $H^\ast({\mathfrak n}(n+1))$ is parametrized by the integral points of the ellipsoid $x_1^2+\dots+x_n^2=x_1x_2+x_2x_3+\dots+x_{n-1}x_n+x_1+\dots+x_n$, or, still better, by the elements of the Weyl group $S_{n+1}$ whose action on the ellipsoid above is generated by the reflections $s_i(x_1,\dots,x_n)=(x_1,\dots,x_{i-1},-x_i+x_{i-1}+x_{i+1}+1,x_{i+1},\dots,x_n),\ i=1,\dots,n$ (in this formula, $x_0$ and $x_{n+1}$ are taken to be zero). If $(p_1,\dots,p_n)=\sigma(0,\dots,0),\, \sigma\in S_{n+1}$, then the corresponding cohomology class $\gamma_\sigma\in H^\ell({\mathfrak n}(n+1))$ where $\ell$ is the length of $\sigma$. In this case, $(p_1,\dots,p_n)$ has a unique presentation as the sum $\sum_{s=1}^q\{i_s,j_s\}$ of different points of the form $$\{i,j\}=(0,\dots,0,\mathop{1}\limits_{(i)},\dots,1,\mathop{0}\limits_{(j)},\dots,0),\ 1\le i<j\le n+1$$ and the class $\gamma_\sigma$ is represented by the monomial cochain $c_{i_1,j_1}\dots c_{i_q,j_q}$ where $c_{i,j}$ takes the value 1 on the one-entry matrix $E_{i,j}$ and takes the value 0 on all other one-entry matrices. Moreover, if the presentations $\sigma(0,\dots,0)=\sum\{i_s,j_s\}, \sigma'(0,\dots,0)=\sum\{i'_t,j'_t\}$ are disjoint and $\sum\{i_s,j_s\}+\sum\{i'_t,j'_t\}=\tau(0,\dots,0)$, then $\gamma_\sigma\gamma_{\sigma'}=\gamma_\tau$; in all other cases, $\gamma_\sigma\gamma_{\sigma'}=0$.

For example, there are 6 permutations in $S_3$: $\sigma_1=(1,2,3),\sigma_2=(2,1,3),\sigma_3=(1,3,2),\sigma_4=(2,3,1),\sigma_5=(3,1,2),\sigma_6=(3,2,1)$. Accordingly, there are 6 integral points on the ellipse $x^2+y^2-x-y-xy=0$,\[\begin{array} {ll} \sigma_1(0,0)=(0,0),&\sigma_4(0,0)=(1,2)=(0,1)+(1,1),\\ \sigma_2(0,0)=(1,0),&\sigma_5(0,0)=(2,1)=(1,0)+(1,1),\\ \sigma_3(0,0)=(0,1),&\sigma_6(0,0)=(2,2)=(1,0)+(0,1)+(1,1), \end{array}\]the cohomology of ${\mathfrak n}(3)$ is spanned by $$\gamma_{\sigma_1}=1\in H^0({\mathfrak n}(3)),\gamma_{\sigma_2},\gamma_{\sigma_3}\in H^1({\mathfrak n}(3)),\gamma_{\sigma_4},\gamma_{\sigma_5}\in H^2({\mathfrak n}(3)),\gamma_{\sigma_6}\in H^3({\mathfrak n}(3)),$$
$\gamma_{\sigma_2}\gamma_{\sigma_4}=-\gamma_{\sigma_3}\gamma_{\sigma_5}=\gamma_{\sigma_6}$, and all other products of cohomology classes of positive dimensions are zero. Similarly for ${\mathfrak n}(4)$ (we write $\sigma_{(ijkl)}$ for the permutation $(i,j,k,l)$): \[\sigma_{(1234)}(0,0,0)=(0,0,0),\ \sigma_{(2134)}(0,0,0)=(1,0,0),\ \sigma_{(1324)}(0,0,0)=(0,1,0),\]
\[ \begin{array} {ll} \sigma_{(1243)}(0,0,0)=(0,0,1),&\sigma_{(3214)}(0,0,0)=(1,0,0)+(0,1,0)+(1,1,0),\\ \sigma_{(2314)}(0,0,0)=(0,1,0)+(1,1,0),&\sigma_{(2341)}(0,0,0)=(0,0,1)+(0,1,1)+(1,1,1),\\ \sigma_{(3124)}(0,0,0)=(1,0,0)+(1,1,0),&\sigma_{(3142)}(0,0,0)=(1,0,0)+(0,0,1)+(1,1,1),\\ \sigma_{(2143)}(0,0,0)=(1,0,0)+(0,0,1),&\sigma_{(2413)}(0,0,0)=(0,1,0)+(1,1,0)+(0,1,1),\\ \sigma_{(1342)}(0,0,0)=(0,0,1)+(0,1,1),&\sigma_{(4123)}(0,0,0)=(1,0,0)+(1,1,0)+(1,1,1),\\ \sigma_{(1423)}(0,0,0)=(0,1,0)+(0,1,1),&\sigma_{(1432)}(0,0,0)=(0,1,0)+(0,0,1)+(0,1,1), \end{array}\] \[\begin{array} {c} \sigma_{(3241)}(0,0,0)=(1,0,0)+(0,0,1)+(0,1,1)+(1,1,1),\\ \sigma_{(2431)}(0,0,0)=(0,1,0)+(0,0,1)+(0,1,1)+(1,1,1),\\ \sigma_{(3412)}(0,0,0)=(0,1,0)+(1,1,0)+(0,1,1)+(1,1,1),\\ \sigma_{(4213)}(0,0,0)=(1,0,0)+(0,1,0)+(1,1,0)+(1,1,1),\\ \sigma_{(4132)}(0,0,0)=(1,0,0)+(0,0,1)+(1,1,0)+(1,1,1),\\  \\ \sigma_{(3421)}(0,0,0)=(0,1,0)+(0,0,1)+(1,1,0)+(0,1,1)+(1,1,1),\\ \sigma_{(4231)}(0,0,0)=(1,0,0)+(0,0,1)+(1,1,0)+(0,1,1)+(1,1,1),\\ \sigma_{(4312)}(0,0,0)=(1,0,0)+(0,1,0)+(1,1,0)+(0,1,1)+(1,1,1),\\ \\ \sigma_{(4321)}(0,0,0)=(1,0,0)+(0,1,0)+(0,0,1)+(1,1,0)+(0,1,1)+(1,1,1).\end{array}\] The cohomology classes of the corresponding monomial cochains form a basis in the cohomology: \[ \begin{array} {c} \gamma_{(1234)}=1\in H^0({\mathfrak n}(4)),\ \gamma_{(2134)},\gamma_{(1324)},\gamma_{(1243)}\in H^1({\mathfrak n}(4)),\\ \gamma_{(2314)},\gamma_{(3124)},\gamma_{(2143)},\gamma_{(1342)},\gamma_{(1423)}\in H^2({\mathfrak n}(4)),\\ \gamma_{(3214)},\gamma_{(2341)},\gamma_{(3142)},\gamma_{(2413)},\gamma_{(4123)},\gamma_{(1432)}\in H^3({\mathfrak n}(4)),\\ \gamma_{(3241)},\gamma_{(2413)},\gamma_{(3412)},\gamma_{(4213)},\gamma_{(4132)}\in H^4({\mathfrak n}(4)),\\ \gamma_{(3421)},\gamma_{(4231)},\gamma_{(4312)}\in H^5({\mathfrak n}(4)),\ \gamma_{(4321)}\in H^6({\mathfrak n}(4)). \end{array} \] The multiplication is described by the following relations: \[\begin{array} {c} \gamma_{(2134)}\gamma_{(1243)}=\gamma_{(2143)};\\ \gamma_{(2134)} \gamma_{(2314)}=-\gamma_{(1324)} \gamma_{(3124)}=\gamma_{(3214)}, \gamma_{(1324)} \gamma_{(1342)}=-\gamma_{(1243)} \gamma_{(1423)}=\gamma_{(1432)};\\ \gamma_{(2134)} \gamma_{(2341)}= \gamma_{(3241)}, \gamma_{(1324)} \gamma_{(2341)}= \gamma_{(2431)},\hskip1in\\ \hskip1in\gamma_{(1324)} \gamma_{(4123)}= -\gamma_{(4213)}, \gamma_{(1243)} \gamma_{(4123)}= -\gamma_{(4132)};\\ \gamma_{(1243)} \gamma_{(3412)}= \gamma_{(2314)} \gamma_{(2341)}= -\gamma_{(3421)},  -\gamma_{(3124)} \gamma_{(2341)}= \gamma_{(1342)} \gamma_{(4123)}= \gamma_{(4231)},\hskip.6in\\ \hskip2.5in  \gamma_{(2134)} \gamma_{(3412)}= \gamma_{(1423)} \gamma_{(4123)}= \gamma_{(4312)};\\ -\gamma_{(2134)} \gamma_{(1243)} \gamma_{(3412)}= -\gamma_{(2134)} \gamma_{(2314)} \gamma_{(2341)}= \gamma_{(1324)} \gamma_{(3124)} \gamma_{(2341)}\hskip.7in\\ \hskip.8in= \gamma_{(1324)} \gamma_{(1342)} \gamma_{(4123)}= -\gamma_{(1243)} \gamma_{(1423)} \gamma_{(4123)}= \gamma_{(3142)} \gamma_{(2413)}= \gamma_{(4321)}. \end{array}\]

Although the procedure given always determines the multiplication in $H^\ast(N(A))$, it does not give a satisfactory explicit description even of the ring $H^\ast({\mathfrak n}(n))$ which remains unclear to us.


\begin{thebibliography}{1}
\bibitem[1]{Fei} Feigin, B.L., Cohomology of current Lie groups and current Lie algebras. (Russian.) {\it Uspekhi mat. nauk}, {\bf 35} (1980), no. 2, 225--226.
\bibitem[2]{Fuchs} Fuchs, D., Cohomology of infinite-dimensional Lie algebras. {\it Consultants Bureau}, NY, 1986. 
\bibitem[3]{GFF} Gelfand, I.M., Feigin, B.L., Fuchs, D.B., Cohomology of infinite-dimensional Lie algebras and Laplace operators. {\it Func. Anal. and Appl.} {\bf 12} (1978), 243--247.
\bibitem[4]{Kac} Kac, V.G., Infinite-dimensional Lie algebras. {\it Cambridge University Press}, Cambridge, 1990.
\bibitem[5]{KK} Kac V.G., Kazhdan D.A., Structure of representations with highest weight of infinite-dimensional Lie algebras. {\it Adv. Math}, {\bf34} (1979), 97--108.
\bibitem[6]{Le} Lepowsky, J., Generalized Verma modules, loop space cohomology, and MacDonald-type identities. {\it Ann. scient. Ec. Norm. Sup.}, {\bf 12} (1979), 169--234.
\bibitem[7]{Wein} Weinstein, F.V., Filtering bases: a tool to compute cohomologies of abstract subalgebras of the Witt algebra. {\it Unconventional Lie algebras}, 155--216. {\it Adv. Soviet Math.}, {\bf 17}, AMS, Providence RI, 1993.
 \end{thebibliography}
\end{document}